\documentclass[11pt]{article}
\usepackage{graphicx}
\usepackage{epsfig}
\usepackage{latexsym}
\newsavebox{\toy}
\savebox{\toy}{\framebox[0.65em]{\rule{0cm}{1ex}}}
\newcommand{\QED}{\usebox{\toy}\end{demo}}
\newenvironment{property}%
{\begin{list}{}{\setlength{\rightmargin}{0pt}%
\setlength{\itemsep}{0pt}}}{\end{list}}
\newlength{\templength}
\newcommand{\bp}{\setlength{\templength}{\labelwidth}%
\setlength{\labelwidth}{2em}\begin{property}}
\newcommand{\ep}{\end{property}\setlength{\labelwidth}{\templength}}
\newtheorem{theorem}{Theorem}[subsection]
\newtheorem{lemma}[theorem]{Lemma}
\newtheorem{proposition}[theorem]{Proposition}
\newtheorem{corollary}[theorem]{Corollary}
\newtheorem{assumption}{Assumption}
\newtheorem{definition}{Definition}[subsection]
\newtheorem{remark}{Remark}[subsection]
\newtheorem{exercise}{Exercise}[subsection]

\newcommand{\Thm}[1]{Theorem \ref{Thm.#1}}
\newcommand{\Lem}[1]{Lemma \ref{Lem.#1}}

\newcommand{\Theorem}[1]{\begin{theorem}\label{Thm.#1}}
\newcommand{\Lemma}[1]{\begin{lemma}\label{Lem.#1}}
\newcommand{\Proposition}[1]{\begin{proposition}\label{Prop.#1}}
\newcommand{\Corollary}[1]{\begin{corollary}\label{Cor.#1}}
\newcommand{\Assumption}[1]{\begin{assumption}\label{Ass.#1}\rm}
\newcommand{\Definition}[1]{\begin{definition}\label{Def.#1}\rm}
\newcommand{\Remark}[1]{\begin{remark}\label{Rem.#1}\rm }
\newcommand{\Exercise}[1]{\begin{exercise}\label{Exe.#1}\rm }

\newcommand{\bd}{\begin{displaymath}}
\newcommand{\ed}{\end{displaymath}}
\newcommand{\bdn}{\begin{equation}}
\newcommand{\bdnl}{\begin{equation}\label}
\newcommand{\edn}{\end{equation}}
\newcommand{\barray}{\begin{array}}
\newcommand{\earray}{\end{array}}
\newcommand{\bds}{\begin{description}}
\newcommand{\eds}{\end{description}}
\newcommand{\bitemize}{\begin{itemize}}
\newcommand{\eitemize}{\end{itemize}}
\newcommand{\benumerate}{\begin{enumerate}}
\newcommand{\eenumerate}{\end{enumerate}}
\newcommand{\btabbing}{\begin{tabbing}}
\newcommand{\etabbing}{\end{tabbing}}
\newcommand{\bcenter}{\begin{center}}
\newcommand{\ecenter}{\end{center}}
\newcommand{\bflushright}{\begin{flushright}}
\newcommand{\bflushleft}{\begin{flushleft}}
\newcommand{\eflushright}{\end{flushright}}
\newcommand{\eflushleft}{\end{flushleft}}
\newcommand{\bdnn }{\begin{eqnarray*}}
\newcommand{\ednn }{\end{eqnarray*}}
\newcommand{\bdmn}{\begin{eqnarray}}
\newcommand{\edmn}{\end{eqnarray}}

\newcommand{\SSC}[1]{\section{#1}\setcounter{equation}{0}}

\newcounter{biblio}
\newenvironment{references}%
{\begin{list}{[\arabic{biblio}]}{\usecounter{biblio}%
\setlength{\leftmargin}{2.5em}\setlength{\rightmargin}{0pt}%
\setlength{\labelwidth}{2em}\setlength{\itemsep}{0pt}}}{\end{list}}
\newcommand{\References}%
{\vspace{2.8ex plus .3ex minus .3ex}%
\begin{center}{\bf References}\end{center}\begin{references}}

\usepackage{amssymb}
\usepackage{mathrsfs}
\newcommand{\bL}{{\mathbb{L}}}
\newcommand{\N}{{\mathbb{N}}}
\newcommand{\Z}{{\mathbb{Z}}}
\newcommand{\zd}{\Z^d}

\newcommand{\R}{{\mathbb{R}}}
\newcommand{\rd}{\R^d}

\newcommand{\ra }{\rightarrow }
\newcommand{\lra }{\longrightarrow }

\newcommand{\Ra}{\Rightarrow }

\newcommand{\LRa}{\Leftrightarrow }
\newcommand{\ov}{\overline}
\newcommand{\tl}{\widetilde}

\newcommand{\Llra}{\Longleftrightarrow }

\newcommand{\vvs}{\vspace{2ex}}
\newcommand{\vs}{\vspace{1ex}}

\newcommand{\lan}{\langle \:}
\newcommand{\ran}{\: \rangle}
\newcommand{\lef}{\left}
\newcommand{\rig}{\right}
\newcommand{\ri}{\right}
\newcommand{\st}{\stackrel}

\newcommand{\8}{\infty}

\newcommand{\dps}{\displaystyle}
\newcommand{\sub}{\subset}

\newcommand{\half}{\mbox{$\frac{1}{2}$}}

\renewcommand{\b}{\beta}

\newcommand{\del}{\delta}

\newcommand{\e}{\varepsilon}

\newcommand{\z}{\zeta}
\newcommand{\h}{\eta}
\newcommand{\tht}{\theta}

\newcommand{\lm}{\lambda}

\newcommand{\n}{\nu}

\newcommand{\s}{\sigma}

\renewcommand{\t}{\tau}

\newcommand{\W}{\Omega}

\newcommand{\cF }{{\cal F}}

\makeatletter
\def\section{\@startsection{section}{1}{\z@}{-3.5ex plus -1ex minus 
 -.2ex}{2.3ex plus .2ex}{\bf}}
\def\subsection{\@startsection{subsection}{2}{\z@}{-3.25ex plus -1ex minus 
 -.2ex}{1.5ex plus .2ex}{\bf}}
\makeatother

\begin{document}
\bcenter

\large{\bf A Note on the Diffusive Scaling Limit 
for a Class of Linear Systems}\footnote{\today}

\vvs \normalsize

\vvs
Yukio Nagahata \\
Department of Mathematics, \\
Graduate School of Engineering Science \\
Osaka University,\\
Toyonaka 560-8531, Japan.\\
email: {\tt nagahata@sigmath.es.osaka-u.ac.jp}\\
URL: {\tt http://www.sigmath.osaka-u.ac.jp/}$\widetilde{}$ {\tt nagahata/}

\vvs
Nobuo Yoshida\footnote{
Supported in part by JSPS Grant-in-Aid for Scientific
Research, Kiban (C) 21540125} \\
Division of Mathematics \\
Graduate School of Science \\
Kyoto University,\\
Kyoto 606-8502, Japan.\\
email: {\tt nobuo@math.kyoto-u.ac.jp}\\
URL: {\tt http://www.math.kyoto-u.ac.jp/}$\widetilde{}$ {\tt nobuo/}

\ecenter
\begin{abstract}
We consider a class of continuous-time stochastic growth models on 
$d$-dimensional lattice with non-negative real numbers 
as possible values per site. 
We remark that the diffusive scaling limit proven in our previous 
work \cite{NY09a} can be extended to wider class of models 
so that it covers the cases 
of potlatch/smoothing processes.
\end{abstract}

\small

\tableofcontents

\normalsize
\SSC{Introduction}
We write $\N^*=\{1,2,...\}$, $\N=\{0\} \cup \N^*$, 
 and 
$\Z=\{ \pm x \; ; \; x \in \N \}$. For 
$x=(x_1,..,x_d) \in \rd$, $|x|$ stands for the $\ell^1$-norm: 
$|x|=\sum_{i=1}^d|x_i|$. For $\h=(\h_x)_{x \in \zd} \in \R^{\zd}$, 
$|\h |=\sum_{x \in \zd}|\h_x|$. 
Let 
$(\W, \cF, P)$ be 
a probability space.
We write $P[X:A]=\int_A X \; dP$  and 
$P[X]=P[ X : \W ]$ for a random variable $X$ and an event $A$. 

\subsection{The model} \label{sec:model}
We go directly into the formal definition of the model, 
referring the reader to \cite{NY09a,NY09b} for relevant backgrounds. 
The class of growth models considered here is a reasonably 
ample subclass of the 
one considered in \cite[Chapter IX]{Lig85} as ``linear systems''. 
We introduce a random vector 
$K=(K_x)_{x \in \zd }$ such that
\bdmn 
& & 0 \le  K_x \le b_K {\bf 1}_{\{ |x| \le r_K\}}\; \; 
\mbox{a.s. for some constants $b_K,r_K \in [0,\8)$,} \label{K_x} \\
& & \mbox{the set $\{x \in \zd\; ;\; P[K_x]\neq 0\}$ contains a 
linear basis of $\rd$.} \label{true_d}
\edmn
The first condition (\ref{K_x}) amounts to the standard boundedness 
and the finite range assumptions for the transition rate of 
interacting particle systems. 
The second condition (\ref{true_d}) makes the model 
``truly $d$-dimensional". 

Let $\t^{z,i}$, ($z \in \zd$, $i \in \N^*$) 
be i.i.d. mean-one exponential random variables 
and $T^{z,i}=\t^{z,1}+...+\t^{z,i}$. 
Let also $K^{z,i}=(K_x^{z,i})_{x \in \zd }$ 
($z \in \zd$, $i \in \N^*$) be i.i.d. 
random vectors with the same distributions as $K$, 
independent of $\{\t^{z,i} \}_{z \in \zd, i \in \N^*}$. 
We suppose that the process $(\h_t)$ 
starts from a deterministic configuration 
$\h_0=(\h_{0,x})_{x \in \zd}\in \N^{\zd}$  
with $|\h_0|<\8$. 
At time $t=T^{z,i}$, $\h_{t-}$ is replaced by $\h_t$, where 
\bdnl{h_(x,t)}
\h_{t,x}=\lef\{ \barray{ll}
K^{z,i}_0\h_{t-,z} & \mbox{if $x=z$}, \\
\h_{t-,x}+K^{z,i}_{x-z}\h_{t-,z} & 
\mbox{if $x \neq z$}.
\earray \ri.
\edn
We also consider the {\it dual process}
$\z_t \in [0,\8)^{\zd}, \; \; t \ge 0$ which evolves in the 
same way as $(\h_t)_{t \ge 0}$ except that 
(\ref{h_(x,t)}) is replaced by its transpose:
\bdnl{h_(x,t)*}
\z_{t,x}=\lef\{ \barray{ll}
\sum_{y \in \zd}K^{z,i}_{y-x}\z_{t-,y} & \mbox{if $x=z$}, \\
\z_{t-,x} & \mbox{if $x \neq z$}.
\earray \ri.
\edn 
Here are some typical examples which fall into the above set-up:

\vs
\noindent $\bullet$
{\bf The binary contact path process (BCPP):}
The binary contact path process (BCPP), originally introduced by 
D. Griffeath \cite{Gri83} is a special case the model, where 
\bdnl{binK} 
K= \lef\{ \barray{ll}
\lef(\del_{x,0}+\del_{x,e} \rig)_{x \in \zd} & 
\mbox{with probability ${\lm \over 2d\lm +1}$, for each $2d$ neighbor 
$e$ of 0} \\
0 & \mbox{with probability ${1 \over 2d\lm +1}$}. 
\earray  \rig. \edn
The process is interpreted as the spread of an infection, 
with $\h_{t,x}$ infected individuals at time $t$ at the site $x$. 
The first line of (\ref{binK}) says that, 
with probability ${\lm \over 2d\lm +1}$ for each $|e|=1$, 
all the infected individuals at site $x-e$ are 
duplicated and added to those on the site $x$.
On the other hand, the second line of (\ref{binK}) says that, 
all the infected individuals at a site become healthy with 
probability ${1 \over 2d\lm +1}$. 
A motivation to study the BCPP comes from the 
fact that the projected process
$\lef( \h_{t,x} \wedge 1\ri)_{x \in \zd},\; \; \; t \ge 0$
is the basic contact process \cite{Gri83}.

\vs
\noindent $\bullet$
{\bf The potlatch/smoothing processes:}
The potlatch process discussed in e.g. 
\cite{HL81} and \cite[Chapter IX]{Lig85} 
is also a special case of the above set-up, 
in which
\bdnl{binP/S} 
K_x=W k_x,\; \; x \in \zd.
\edn
Here, $k=(k_x)_{x \in \zd} \in [0,\8)^{\zd}$ is a 
non-random vector and $W$ is a non-negative, bounded, mean-one 
random variable such that $P(W = 1)<1$ 
(so that the notation $k$ here is consistent with the definition 
(\ref{kp12}) below). The smoothing process is the dual process 
of the potlatch process. 
The potlatch/smoothing processes 
were first introduced in 
\cite{Spi81} for the case $W \equiv 1$ and discussed further in 
\cite{LS81}. It was in \cite{HL81} where case with $W \not\equiv 1$ was 
introduced and discussed. 
Note that we {\it do not} assume that $k_x$ is a 
transition probability of an irreducible random walk, 
unlike in the literatures mentioned above. 

\vvs
We now recall the following facts from 
\cite[page 433, Theorems 2.2 and 2.3]{Lig85}.
  Let $\cF_t$ be the 
$\s$-field generated by $\h_s$, $s \le t$.
Let $(\h^x_t)_{t \ge 0}$ be the process 
$(\h_t)_{t \ge 0}$ starting from one particle 
at the site $x$: $\h^x_0=\del_x$.
Similarly, let 
$(\z^x_t)_{t \ge 0}$ be the dual process starting from one particle 
at the site $x$: $\z^x_0=\del_x$.
\Lemma{0,1}
We set:
\bdmn
k&=& (k_x)_{x \in \zd}=(P[K_x])_{x \in \zd} \; \; \; \label{kp12} \\
\ov{\h}_t &=&(e^{-(|k|-1)t }\h_{t,x})_{x \in \zd}. \label{ovh_t}
\edmn
Then,
\bds
\item[a)]
$(|\ov{\h}_t|, \cF_t)_{t \ge 0}$ is a martingale, 
and therefore, the following limit exists a.s. 
\bdnl{ovn_8}
|\ov{\h}_\8| =\lim_{t \ra \8}|\ov{\h}_t|.
\edn
\item[b)] Either 
\bdnl{0,1}
P[|\ov{\h}^0_\8|]=1\; \; \mbox{or}\; \; 0.
\edn
Moreover, $P[|\ov{\h}^0_\8|]=1$ if and only if the limit (\ref{ovn_8}) 
is convergent in $\bL^1 (P)$. 
\item[c)] The above a)--b), with $\h_t$ replaced by $\z_t$ 
are true for the dual process.
\eds
\end{lemma}
\subsection{Results}
We are now in position to state our main result in this article (\Thm{CLT}). 
It extends our previous result \cite[Theorem 1.2.1]{NY09a} 
to wider class of models 
so that it covers the cases 
of potlatch/smoothing processes, cf. Remarks 1)--2) after \Thm{CLT}.

We first introduce some more notation. 
For $\h, \z \in \R^{\zd}$, the inner product and the 
discrete convolution are  defined respectively by 
\bdnl{<>*}
\lan \h, \z \ran =\sum_{x \in \zd}\h_x\z_x
\; \; \; \mbox{and}\; \; \;
(\h * \z)_x  =\sum_{y \in \zd}\h_{x-y}\z_y
\edn
provided the summations converge. 
We define for 
$x,y \in \zd$,
\bdnl{beta_x}
\b_{x,y}=P[(K-\del_0)_x(K-\del_0)_y]\;  \; \; \mbox{and}\; \; \; 
\b_x = \sum_{y \in \zd}\b_{x+y,y}
\edn
If we simply write $\b$ in the sequel, 
it stands for the function $x \mapsto \b_x$.
Note then that
\bdnl{|beta|}
\lan \b, 1 \ran = \sum_{x,y \in \zd}\b_{x,y}= P[(|K|-1)^2].
\edn
We also introduce:
\bdnl{G(x)}
G_S(x)=\int^\8_0P_S^0(S_t=x)dt,
\edn
where $((S_t)_{t \ge 0}, P_S^x)$ is the continuous-time random walk on 
$\zd$ starting from $x \in \zd$, with the generator 
\bdnl{L_S}
L_Sf (x)=\sum_{y \in \zd}L_S (x,y)
\lef( f(y)-f(x)\rig), \; \; \mbox{with} \; \; 
L_S (x,y)={ k_{x-y}+k_{y-x} \over 2}
\; \; \mbox{for $x \neq y$},
\edn
cf. (\ref{kp12}).
The set of bounded continuous functions on $\rd$ is 
denoted by $C_{\rm b} (\rd)$.
\Theorem{CLT}
Suppose $d \ge 3$. 
Then, the following conditions are equivalent:
\bds 
\item[a)] 
$\lan \b, G_S \ran <2.$
\item[b)] There exists a bounded 
function $h:\zd \ra [1,\8)$ 
such that:
\bdnl{eq:h}
(L_S h)(x)+\half \del_{0,x}\lan \b ,h \ran  \le 0, \; \; \; x\in \zd.
\edn
\item[c)] 
${\dps \sup_{t \ge 0} P[|\ov{\h}_t|^2]<\8}$.
\item[d)]
${\dps  \lim_{t \ra \8} \sum_{x \in \zd}
f\lef((x -mt)/\sqrt{t}\ri)\ov{\h}_{t,x}=|\ov{\h}_\8|\int_{\rd}fd\n}$  
in $\bL^2 (P)$ for all $f \in C_{\rm b} (\rd)$, \\
where $ m=\sum_{x \in \zd}xk_x \in \rd $ 
and $\n$ is the Gaussian measure with
\bdnl{nu}
\int_{\rd}x_id\n (x)=0, 
\; \; \; \int_{\rd}x_ix_jd\n (x)
=\sum_{x \in \zd}x_ix_jk_x,\; \; \; 
i,j=1,..,d.
\edn
\item[b')] There exists a bounded 
function $h:\zd \ra [1,\8)$ 
such that:
\bdnl{eq:h*}
(L_S h)(x)+\half h (0) \b_x \le 0, \; \; \; x\in \zd.
\edn
\item[c')] 
${\dps \sup_{t \ge 0} P[|\ov{\z}_t|^2]<\8}$.
\item[d')]
${\dps  \lim_{t \ra \8} \sum_{x \in \zd}
f\lef((x -mt)/\sqrt{t}\ri)\ov{\z}_{t,x}=|\ov{\z}_\8|\int_{\rd}fd\n}$  
in $\bL^2 (P)$ for all $f \in C_{\rm b} (\rd)$.
\eds
\end{theorem}
The main point of \Thm{CLT} is that a) implies d) and d'), 
while the equivalences between the other conditions 
are byproducts. 

\vvs
\noindent {\bf Remarks}:
\noindent {\bf 1)}
\Thm{CLT} extends \cite[Theorem 1.2.1]{NY09a}, where 
the following extra technical condition was imposed:
\bdnl{assum_09a}
\b_x=0 \; \; \mbox{for $x \neq 0$.}
\edn
For example, BCPP satisfies (\ref{assum_09a}), 
while the potlatch/smoothing processes do not. \\
\noindent {\bf 2)}
Let $\pi_d$ be the return probability for the simple random 
walk on $\zd$.
We then have that
\bdnl{cond:sloc_exam}
\lan \b, G_S \ran <2\; \; \Llra \; \; 
\lef\{ \barray{ll}
\lm >{1 \over 2d (1-2\pi_d)} & \mbox{for BCPP},\\
P[W^2] <{(2|k|-1)G_S(0)\over \lan G_S*k, k \ran}
& \mbox{for the potlatch/smoothing processes.} 
\earray \rig.
\edn
cf. \cite[page 460, (6.5) and page 464, Theorem 6.16 (a)]{Lig85}.
For BCPP, (\ref{cond:sloc_exam}) can be seen from that 
(cf. \cite[page 965]{NY09a})
$$
\b_{x,y}={{\bf 1}\{ x=0\}+\lm {\bf 1}\{ |x|=1\} \over 2d\lm +1}\del_{x,y}, 
\; \; \mbox{and}\; \; 
G_S(0)={2d\lm +1 \over 2d\lm}{1 \over 1-\pi_d}.
$$
To see (\ref{cond:sloc_exam}) for the potlatch/smoothing processes, 
we note that $\half (k+ \check{k})*G_S =|k|G_S-\del_0$, 
with $\check{k}_x=k_{-x}$ and that
$$
\b_{x,y}=P[W^2]k_xk_y-k_x \del_{y,0}-k_y \del_{x,0}+\del_{x,0}\del_{y,0}.
$$ 
Thus, 
\bdnn
\lan \b, G_S \ran
& = & P[W^2]\lan G_S*k, k \ran -\lan G_S, k+ \check{k} \ran +G_S(0) \\
& = & P[W^2]\lan G_S*k, k \ran +2-(2|k|-1)G_S(0),
\ednn
from which (\ref{cond:sloc_exam}) 
for the potlatch/smoothing processes follows. \\
\noindent {\bf 3)}
It will be seen from the proof that 
the inequalities in (\ref{eq:h}) and (\ref{eq:h*}) can be replaced by 
the equality, keeping the other statement of \Thm{CLT}.

\vvs
As an immediate consequence of \Thm{CLT}, we have the following 
\Corollary{CLT}
Suppose either of a)--d) in \Thm{CLT}.
Then, $P [|\ov{\h}_\8|]=|\h_0|$ and for all $f \in C_{\rm b} (\rd)$, 
\bdnn 
& & \lim_{t \ra \8} \sum_{x \in \zd}
f\lef((x -mt)/ \sqrt{t}\ri)\frac{\h_{t,x}}{|\h_t|}
{\bf 1}_{\{ \h_t \not\equiv 0 \}}
=\int_{\rd}fd\n \\
& & 
\mbox{in probability with respect to 
$P(\; \cdot \; | \h_t \not\equiv 0, \; \forall t)$}.
\ednn
where $ m=\sum_{x \in \zd}xk_x \in \rd $ and 
$\n$ is the Gaussian measure defined by (\ref{nu}). 
Similarly, either of a),b'),c'),d') in \Thm{CLT} 
implies the above statement, 
with  $\h_t$ replaced by the dual process $\z_t$.
\end{corollary}
Proof: The case of $(\h_\cdot )$ follows from \Thm{CLT}d). Note also that 
if $P (|\ov{\h}_\8| >0)>0$, then, up to a null set, 
$$
\{ \; |\ov{\h}_\8| >0 \; \}=\{\; \h_t \not\equiv 0, \; \forall t\; \},
$$
which follows from 
\cite[Lemma 2.1.2]{NY09b}. 
The proof for the case of $(\z_\cdot )$ is the same.
\hfill $\Box$
\SSC{The proof of \Thm{CLT}}
\subsection{The equivalence of a)--c)} \label{sec:(a)--(c)}
We first show the Feynman-Kac formula for two-point 
function, which is the basis of the proof of \Thm{CLT}. 
To state it, we introduce Markov chains 
$(X,\tl{X})$ and $(Y,\tl{Y})$ which are also exploited in \cite{NY09a}. 
Let 
$(X,\tl{X})=((X_t, \tl{X}_t)_{t \ge 0}, P_{X,\tl{X}}^{x,\tl{x}})$ and 
$(Y,\tl{Y})=((Y_t, \tl{Y}_t)_{t \ge 0}, P_{Y,\tl{Y}}^{x,\tl{x}})$ be 
the continuous-time Markov chains on 
$\zd \times \zd$ starting from $(x,\tl{x})$, with the generators  
\bdnl{LYY}
\barray{ll}
 & {\dps L_{X,\tl{X}}f (x,\tl{x})=\sum_{y,\tl{y} \in \zd}
L_{X,\tl{X}}(x,\tl{x}, y,\tl{y})  \lef( f(y,\tl{y})-f(x,\tl{x})\rig)},\\
\mbox{and} & {\dps L_{Y,\tl{Y}}f (x,\tl{x})=\sum_{y,\tl{y} \in \zd}
L_{Y,\tl{Y}}(x,\tl{x}, y,\tl{y})  \lef( f(y,\tl{y})-f(x,\tl{x})\rig)},
\earray 
\edn
respectively, where
\bdnl{Gm2}
\barray{ll}
 & {\dps
L_{X,\tl{X}}(x,\tl{x}, y,\tl{y})
 =(k-\del_0)_{x-y}\del_{\tl{x},\tl{y}}
+(k-\del_0)_{\tl{x}-\tl{y}}\del_{x,y}
+\b_{x-y, \tl{x}-y}\del_{y,\tl{y}} }\\
\mbox{and} & L_{Y,\tl{Y}}(x,\tl{x}, y,\tl{y})=L_{X,\tl{X}}(y,\tl{y},x,\tl{x})
\earray \edn 
It is useful to note that 
\bdmn
\sum_{y,\tl{y}}L_{X,\tl{X}}(x,\tl{x}, y,\tl{y})
&= & 2(|k|-1)+\b_{x-\tl{x}}, \label{sumL} \\
\sum_{y,\tl{y}}L_{Y,\tl{Y}}(x,\tl{x}, y,\tl{y})
&= & 2(|k|-1)+\lan \b, 1 \ran  \del_{x,\tl{x}}.\label{sumL*} 
\edmn
Recall also the notation $(\h^x_t)_{t \ge 0}$ and $(\z^x_t)_{t \ge 0}$ 
introduced before \Lem{0,1}.
\Lemma{FK2}
For $t \ge 0$ and $x,\tl{x}, y,\tl{y}\in \zd$, 
\bdnn
P[\z^y_{t,x}\z^{\tl{y}}_{t,\tl{x}}]
& = &  P[\h^x_{t,y}\h^{\tl{x}}_{t,\tl{y}}]  \\ 
& = &  
e^{2(|k|-1)t}P^{y,\tl{y}}_{X,\tl{X}}
\lef[ e_{X,\tl{X},t}: 
(X_t, \tl{X}_t)=(x,\tl{x})\;  \ri]  \\ 
& = & e^{2(|k|-1)t}P^{x,\tl{x}}_{Y,\tl{Y}}
\lef[ e_{Y,\tl{Y},t}: 
(Y_t,\tl{Y}_t)=(y,\tl{y})\;  \ri], 
\ednn
where 
$e_{X,\tl{X},t}
=\exp \lef( \int^t_0 \b_{X_s-\tl{X}_s}ds\ri)$ and 
$e_{Y,\tl{Y},t}
=\exp \lef( \lan \b, 1 \ran \int^t_0 \del_{Y_s,\tl{Y}_s}ds\ri)$.
\end{lemma}
Proof:  
By the time-reversal argument as in \cite[Theorem 1.25]{Lig85}, we see that 
$(\h^x_{t,y}, \h^{\tl{x}}_{t,\tl{y}})$ 
and $(\z^y_{t,x}, \z^{\tl{y}}_{t,\tl{x}})$ have the same law.
This implies the first equality.
In \cite[Lemma 2.1.1]{NY09a}, we showed the second equality, 
using (\ref{sumL*}). Finally, we see from 
(\ref{Gm2}) -- (\ref{sumL*}) that the operators:
\bdnn
f (x,\tl{x}) 
& \mapsto & 
L_{X,\tl{X}}f(x,\tl{x})+\b_{x-\tl{x}}f(x,\tl{x}), \\
f (x,\tl{x}) 
& \mapsto & 
L_{Y,\tl{Y}}f(x,\tl{x})+\lan \b, 1 \ran \del_{x,\tl{x}}f(x,\tl{x})
\ednn
are transpose to each other, and hence are the semi-groups generated by the 
above operators. This proves the last equality of the lemma.
\hfill $\Box$

\vvs 
\Lemma{Y-Y}
$((X_t-\tl{X}_t)_{t \ge 0}, P_{X,\tl{X}}^{x,0})$ and 
$((Y_t-\tl{Y}_t)_{t \ge 0}, P_{Y,\tl{Y}}^{x,0})$ are  
Markov chains with the generators:
\bdnl{Y-Y}
\barray{ll}
 & L_{X-\tl{X}}f(x)=2L_S f(x)+\b_x(f(0)-f(x)) \\
\mbox{and} & 
L_{Y-\tl{Y}}f(x)=2L_S f(x)
+(\lan \b, f  \ran -\lan \b, 1  \ran f(x))\del_{x,0},
\earray 
\edn
respectively (cf. (\ref{L_S})). 
Moreover, these Markov chains are transient for 
$d \ge 3$. 
\end{lemma}
Proof: Let $(Z,\tl{Z})=(X,\tl{X})$ or $(Y,\tl{Y})$. 
Since $(Z,\tl{Z})$ is shift-invariant, in the sense that 
$L_{Z,\tl{Z}} (x+v,\tl{x}+v,y+v,\tl{y}+v) 
=L_{Z,\tl{Z}} (x,\tl{x},y,\tl{y}) $ 
for all $v \in \zd$, $((Z_t-\tl{Z}_t)_{t \ge 0}, P_{Z,\tl{Z}}^{x,\tl{x}})$ 
is a Markov chain. Moreover, the jump rates 
$L_{Z-\tl{Z}}(x,y)$, $x \neq y$ are  computed as follows:
$$
L_{Z-\tl{Z}}(x,y) = \sum_{z \in \zd}L_{Z,\tl{Z}} (x,0,z+y,z)
 = \lef\{\barray{ll}
k_{x-y}+k_{y-x} +\del_{y,0}\b_x & \mbox{if  $(Z,\tl{Z})=(X,\tl{X})$,}\\
k_{x-y}+k_{y-x} +\del_{x,0}\b_y & \mbox{if  $(Z,\tl{Z})=(Y,\tl{Y})$.}
\earray \rig.
$$
These prove (\ref{Y-Y}). By (\ref{true_d}), 
the random walk $S_\cdot$ is transient for $d \ge 3$. Thus, 
$Z-\tl{Z}$ is transient $d \ge 3$, since 
$L_{Z-\tl{Z}} (x, \cdot)=2L_S (x, \cdot)$ except for finitely many 
$x$.
\hfill $\Box$

\vvs
\noindent {\it Proof of a) $\LRa$ b) $\LRa$ c)}:
a) $\Ra$ b):
Under the assumption a), 
the function $h$ given below satisfies conditions in b): 
\bdnl{h=1+cG}
h=1+cG_S\; \; \mbox{with}\; \; 
c={ \lan \b, 1 \ran \over 2- \lan \b, G_S \ran }.
\edn
In particular it solves (\ref{eq:h}) with equality. \\
b) $\Ra$ c): By \Lem{FK2}, we have that
\bds \item[1)] \hspace{1cm}
$P[|\ov{\h}^x_t||\ov{\h}^{\tl{x}}_t|]
=P^{x,{\tl{x}}}_{Y,\tl{Y}}
\lef[ e_{Y,\tl{Y},t}  \ri],
\; \; \; x,\tl{x} \in \zd$,
\eds
where 
$e_{Y,\tl{Y},t}=\exp \lef( \lan \b, 1 \ran \int^t_0 \del_{Y_s,\tl{Y}_s}ds\ri)$.
By \Lem{Y-Y}, (\ref{eq:h}) reads:
\bdnl{eq:h2}
L_{Y-\tl{Y}}h(x)+\lan \b, 1 \ran \del_{x,0}h(x) \le 0,\; \; \; x \in \zd
\edn
and thus,
$$
P^{x,0}_{Y,\tl{Y}}
\lef[ e_{Y,\tl{Y},t} h(Y_t-\tl{Y}_t) \ri] \le h(x),
\; \; \; x \in \zd.
$$
Since $h$ takes its values in $[1,\sup h]$ with $\sup h <\8$, we have 
$$
\sup_xP^{x,0}_{Y,\tl{Y}}\lef[ e_{Y,\tl{Y},t} \ri]
\le \sup h <\8.
$$
By this and 1), we obtain that
$$
\sup_x P[|\ov{\h}^x_t||\ov{\h}^0_t|] \le \sup h <\8.
$$
c) $\Ra$ a) : Let $G_{Y-\tl{Y}}(x,y)$ be the Green function 
of the Markov chain $Y-\tl{Y}$ (cf. \Lem{Y-Y}). Then, it follows from 
 (\ref{Y-Y}) that
\bdnl{GG}
G_{Y-\tl{Y}}(x,y)=\half G_S (y-x)
+\half \lef( \lan \b, G_S \ran -\lan \b, 1 \ran G_S (0) \ri)G_{Y-\tl{Y}}(x,0).
\edn
On the other hand, 
we have by 1) that for any $x,\tl{x} \in \zd$,
$$
P^{x,\tl{x}}_{Y,\tl{Y}}
\lef[ e_{Y,\tl{Y},t} \ri]
=P[|\ov{\h}^x_t||\ov{\h}^{\tl{x}}_t|] \le P[|\ov{\h}^0_t|^2],
$$
where the last inequality comes from Schwarz inequality 
and the shift-invariance.  Thus, 
\bdnl{supFK2}
P^{x,{\tl{x}}}_{Y,\tl{Y}}
\lef[ e_{Y,\tl{Y},\8}  \ri]
\le \sup_{t \ge 0}P[|\ov{\h}^0_t|^2]<\8.
\edn
Therefore, we can define $h: \zd \ra [1,\8)$ by:
\bdnl{h(X)=P^x[e_8]}
h(x)=P^{x,0}_{Y,\tl{Y}}
\lef[ e_{Y,\tl{Y},\8}  \ri],
\edn
which solves:
$$
h(x)=1+G_{Y-\tl{Y}}(x,0)\lan \b, 1 \ran h (0).
$$
For $x=0$, it implies that
$$
G_{Y-\tl{Y}}(0,0)\lan \b, 1 \ran <1.
$$
Plugging this into (\ref{GG}), we have a).  
\hfill $\Box$

\vvs
\noindent {\bf Remark:}
The function $h$ defined 
by (\ref{h(X)=P^x[e_8]}) solves (\ref{eq:h2}) with equality, as can be 
seen by the way it is defined. 
This proves c) $\Ra$ b) directly.
It is also easy to see from (\ref{GG}) that the function $h$ defined 
by (\ref{h(X)=P^x[e_8]}) and by (\ref{h=1+cG}) coincide.
\subsection{The equivalence of c) and d)}
To proceed from c) to the diffusive scaling limit d), 
we will use the following variant of \cite[Lemma 2.2.2]{NY09a}:
\Lemma{perturb2}
Let $((Z_t)_{t \ge 0}, P^x)$ be a continuous-time random walk 
on $\zd$ starting from $x$, 
with the generator:
$$
L_Z f(x)=\sum_{y \in \zd}L_Z(x,y)(f(y)-f(x)),
$$ 
where we assume that:
$$
\sum_{x \in \zd}|x|^2L_Z(0,x)<\8.
$$
On the other hand, let $\tl{Z}=((\tl{Z}_t)_{t \ge 0}, \tl{P}^x)$ 
be the continuous-time Markov chain on $\zd$ starting from $x$, 
with the generator:
$$
L_{\tl{Z}} f(x)=\sum_{y \in \zd}L_{\tl{Z}}(x,y)(f(y)-f(x)).
$$ 
We assume that $z \in \zd$, 
$D \sub \zd$ and a function $v:\zd \ra \R$ satisfy
\bdnn
& & \mbox{$L_Z(x,y)=L_{\tl{Z}}(x,y)$ if $x \not\in D \cup \{y\}$,} \\
& & \mbox{$D$ is transient for both $Z$ and $\tl{Z}$,} \\
& & \mbox{$v$ is bounded and $v \equiv 0$ outside $D$,}\\
& & \mbox{$e_t\st{\rm def}{=}\exp \lef( \int_0^t v(\tl{Z}_u)du\ri)$, $t \ge 0$ 
are uniformly integrable with respect to $\tl{P}^z$.}
\ednn
Then, for $f \in C_{\rm b}(\rd)$, 
$$
\lim_{t \ra \8}\tl{P}^z
\lef[e_tf((\tl{Z}_t-mt)/\sqrt{t})\ri]
=\tl{P}^z\lef[e_\8\ri]\int_{\rd}fd\n,
$$
where $m=\sum_{x \in \zd}xL_Z (0,x)$ and $\n$ is the Gaussian measure with:
$$
\int_{\rd}x_id\n (x)=0, 
\; \; \; \int_{\rd}x_ix_jd\n (x)
=\sum_{x \in \zd}x_ix_jL_Z (0,x),\; \; \; 
i,j=1,..,d.
$$
\end{lemma}
Proof: We refer the reader to the 
proof of \cite[Lemma 2.2.2]{NY09a}, which  works 
almost verbatim here. 
The uniform integrability of $e_t$ is used to make sure 
that $\lim_{s \ra \8}\sup_{t \ge 0}|\e_{s,t}|=0$, where $\e_{s,t}$ 
is an error term introduced in the proof of \cite[Lemma 2.2.2]{NY09a}.
\hfill $\Box$

\vvs
\noindent {\it Proof of c) $\LRa$ d):}
c) $\Ra$ d): 
Once (\ref{supFK2}) is obtained, 
we can conclude d) exactly in the same way as 
in the corresponding part of \cite[Theorem 1.2.1]{NY09a}.
Since c) implies that 
$\lim_{t \ra \8}|\ov{\h}_t|=|\ov{\h}_\8|$ in $\bL^2 (P)$,
 it is enough to prove that
$$
U_t\st{\rm def.}{=}\sum_{x \in \zd}\ov{\h}_{t,x}f \lef( (x-mt)/\sqrt{t}\rig)
\lra 0\; \; \; \mbox{in $\bL^2 (P)$ as $t \nearrow \8$}
$$
for $f \in C_{\rm b} (\rd)$ such that  $\int_{\rd}fd\n =0$.
We set $f_t(x, \tl{x}) =f ((x-m)/\sqrt{t})f ((\tl{x}-m)/\sqrt{t})$. 
By \Lem{FK2},
$$
P[U_t^2]
 =  \sum_{x,\tl{x} \in \zd}
P[\ov{\h}_{t,x}\ov{\h}_{t,\tl{x}}]f_t(x, \tl{x})
= \sum_{x,\tl{x} \in \zd}\h_{0,x}\h_{0,\tl{x}}
P_{Y,\tl{Y}}^{x,\tl{x}}\lef[ e_{Y,\tl{Y},t}f_t(Y_t, \tl{Y}_t) \ri].
$$
Note that by (\ref{supFK2}) and c), 
\bds
\item[1)]\hspace{1cm}
${\dps P_{Y,\tl{Y}}^{x,\tl{x}}\lef[ e_{Y,\tl{Y},\8}\ri] <\8}$.
\eds
Since $|\h_0|<\8$, it is enough to prove that for each $x,\tl{x} \in \zd$
$$
\lim_{t \ra \8}
P_{Y,\tl{Y}}^{x,\tl{x}}\lef[ e_{Y,\tl{Y},t}f_t(Y_t, \tl{Y}_t) \ri]=0.
$$
To prove this, we apply \Lem{perturb2} to the Markov chain 
$\tl{Z}_t\st{\rm def.}{=}(Y_t,\tl{Y}_t)$ and the random walk 
$(Z_t)$ on $\zd \times \zd$ with the generator:
$$
L_Zf (x,\tl{x})
=\sum_{y,\tl{y} \in \zd}L_Z (x,\tl{x},y,\tl{y}) 
\lef( f(y,\tl{y})-f(x,\tl{x})\rig),
$$
where 
$$
L_Z (x,\tl{x},y,\tl{y}) =\lef\{\barray{ll}
k_{\tl{y}-\tl{x}}
& \mbox{if $x=y$ and $\tl{x} \neq \tl{y}$,}\\
k_{y-x} 
& \mbox{if $x \neq y$ and $\tl{x} =\tl{y}$,}\\
0 & \mbox{if otherwise}.\earray \rig.
$$
Let $D=\{ (x,\tl{x}) \in \zd \times \zd\; ; \; x=\tl{x}\}$. 
Then, 
\bds
\item[2)] $L_Z (x,\tl{x},y,\tl{y})=L_{Y,\tl{Y}}(x,\tl{x}, y,\tl{y})$ if 
$ (x,\tl{x}) \not\in D \cup \{ (y,\tl{y}) \}$.
\eds
Moreover, by \Lem{Y-Y}, 
\bds
\item[3)] $D$ is transient both for $(Z_t)$ and for $(\tl{Z}_t)$.
\eds
 Finally, the Gaussian measure $\n \otimes \n $ 
is the limit law in the central limit theorem 
for the random walk $(Z_t)$.
Therefore, by 1)--3) and \Lem{perturb2},
$$
\lim_{t \ra \8}
P_{Y,\tl{Y}}^{x,\tl{x}}\lef[ e_{Y,\tl{Y},t}f_t(Y_t, \tl{Y}_t)  \ri]
=P_{Y,\tl{Y}}^{x,\tl{x}}\lef[ e_{Y,\tl{Y},\8}\ri]
\lef( \int_{\rd}fd\n \rig)^2=0.
$$
d) $\Ra$ c):This can be seen by taking $f \equiv 1$. 
\hfill $\Box$ 
\subsection{The equivalence of a),b'),c')}
a) $\Ra$ b'): 
Let 
$h=2- \lan \b, G_S\ran +\b *G_S.$
Then, it is easy to see that $h$ solves 
(\ref{eq:h*}) with equality.
Moreover, using \Lem{b*G} below, we see that $h(x)>0$ for $x \neq 0$ by 
as follows:
\bdnn
(\b *G_S)(x)-(\b *G_S)(0) 
& \ge & 
\lef( {G_S(x)\over G_S(0)}-1 \ri)(\b *G_S)(0) -2 {G_S(x)\over G_S(0)} \\
& > & \lef( {G_S(x)\over G_S(0)}-1 \ri)2 
-2 {G_S(x)\over G_S(0)}=-2.
\ednn
Since $h(0)=2$ 
and $\lim_{|x| \ra \8}h(x)=2-(\b *G_S)(0)  \in (0,\8)$, 
$h$ is bounded away from both $0$ and $\8$. 
Therefore, a constant multiple of the above $h$ satisfies the 
conditions in b'). \\
b') $\LRa$ c'): This can be seen similarly as b) $\LRa$ c) 
(cf. the remark at the end of section \ref{sec:(a)--(c)}). \\
c') $\Ra$ a) : 
We first note that 
\bds \item[1)] \hspace{1cm}
${\dps \lim_{|x| \ra \8} (\beta * G_S)(x)=0 }$,
\eds
since $G_S$ vanishes at infinity and $\b$ is of finite support. 
We then set:
$$
h_0(x) =   P^{x,0}_{X,\tilde{X}}[e_{X,\tilde{X},\infty}],
\; \; \; 
h_2(x)  =  h_0(x) - \frac{1}{2} h_0(0) (\beta * G_S)(x).
$$ 
Then, there exists positive constant $M$ such that
$\frac{1}{M} \le h_0 \le M$ and
$$
(L_S h_0)(x) = - \frac{1}{2} h_0(0) \beta_x,
~~~\textrm{for all } x \in \mathbf{Z}^d.
$$
By 1), $h_2$ is also bounded and 
\[
(L_S h_2)(x) = (L_S h_0)(x) - \frac{1}{2} h_0(0) L_S (\beta * G_S)(x)
= - \frac{1}{2} h_0(0) \beta_x + \frac{1}{2} h_0(0) \beta_x = 0.
\]
This implies that there exists a constant $c$ such that 
$h_2 \equiv c$ on the subgroup $H$ of $\zd$ generated by the 
set $\{ x \in \zd \; ; \; k_x+k_{-x} >0\}$, i.e., 
\bds \item[2)] \hspace{1cm}
${\dps
h_0(x) - \frac{1}{2} h_0(0) (\beta * G_S)(x) = c\; \; }$ for $x \in H$.
\eds 
By setting $x=0$ in 2), we have
\[
c= h_0(0) (1 - \frac{\langle \beta, G_S \rangle}{2}).
\]
On the other hand, we see from 1)--2) that 
\[
0<\frac{1}{M} \le \lim_{|x| \ra \8 \atop x \in H}h_0(x)= c.
\]
These imply $\langle \beta, G_S \rangle < 2$.
\hfill $\Box$
\Lemma{b*G}
For $d \ge 3$,
$$
 (\b *G_S)(x) \ge {G_S(x)\over G_S(0)}\lef( (\b *G_S)(0) -2 \ri)+2\del_{0,x} 
\; \; \; x \in \zd.
$$
\end{lemma}
Proof: The function $\b_x$ can be 
either positive or negative. To control this inconvenience, 
we introduce: $\tl{\b}_x = \sum_{y \in \zd}P[K_yK_{x+y}]$.
Since $\tl{\b}_x \ge 0$ and $G_S(x+y)G_S(0) \ge G_S(x)G_S(y)$ 
for all $x,y \in \zd$, we have 
\bds \item[1)] \hspace{1cm} 
$G_S(0)(G_S * \tl{\b})(x) \ge G_S(x)(G_S * \tl{\b})(0)$.
\eds
On the other hand, it is easy to see that
$$
\b=\tl{\b}-k-\check{k}+\del_0, \; \; \mbox{with}\; \; \check{k}_x=k_{-x}.
$$
Therefore, using $\half (k+ \check{k})*G_S =|k|G_S-\del_0$, 
\bds \item[2)] \hspace{1cm} 
$\b*G_S  =  (\tl{\b}-k-\check{k}+\del_0)*G_S 
=\tl{\b}*G_S -(2|k|-1)G_S +2\del_0$.
\eds
Now, by 1)--2) for $x=0$, 
$$
(G_S * \tl{\b})(x) 
\ge {G_S(x) \over G_S(0)}(G_S * \tl{\b})(0)
={G_S(x) \over G_S(0)}(\b*G_S (0)-2)+(2|k|-1)G_S(x).
$$
Plugging this in 2), we get the desired inequality.
\hfill $\Box$
\subsection{The equivalence of c') and d')}
d') $\Ra$ c'): This can be seen by taking $f \equiv 1$. \\
c') $\Ra$ d'):By \Lem{FK2}, Schwarz inequality and the shift-invariance, 
we have that
$$
P^{x,\tl{x}}_{X,\tl{X}}
\lef[ e_{X,\tl{X},t}  \ri]=P[|\ov{\z}^x_t||\ov{\z}^{\tl{z}}_t|]
 \le P[|\ov{\z}^0_t|^2],\; \; \; \mbox{for $x,\tl{x} \in \zd$},
$$
where $e_{X,\tl{X},t}=\exp \lef( \int^t_0 \b_{X_s-\tl{X}_s}ds\ri)$.
Thus, under c'), the following function is well-defined:
$$
h_0(x)\st{\rm def}{=} P^{x,0}_{X,\tilde{X}}[e_{X,\tilde{X},\infty}].
$$
Moreover, there exists $M \in (0,\8)$ such that
$\frac{1}{M} \le h_0 \le M$ and
$$
(L_S h_0)(x) =- \frac{1}{2} h_0(0) \beta_x,
~~~\textrm{for all } x \in \zd.
$$
We set
$$
h_1 (x) = h_0(x) - \frac{1}{2M}.
$$
Then, we have $0 < \frac{1}{2M} \le h_0 \le M$ and
$$
L_S h_1 (x) = L_S h_0 (x) = - \frac{1}{2} h_0(0) \beta_x
= - \frac{1}{2} h_1 (0) p\beta_x, 
\; \; \; \mbox{with}\; \; \; p={h_0(0) \over h_1(0)}>1.
$$
This implies, as in the proof of b) $\Ra$ c) that 
$$
\sup_{t \ge 0}P^{x,\tl{x}}_{X,\tl{X}}
\lef[ e_{X,\tl{X},t}^p  \ri] \le 2M^2 <\8
\; \; \; \mbox{for $x,\tl{x} \in \zd$},
$$
which guarantees the uniform integrability of 
$e_{X,\tl{X},t}$, $t \ge 0$ required to apply \Lem{perturb2}.
The rest of the proof is the same as in c) $\Ra$ d).
\hfill $\Box$

\small


\end{document}